\documentclass[12pt]{article}
\usepackage[amsmath]{e-jc}

\usepackage{tabularx}
\usepackage{graphicx}

\dateline{-----}{--------}{TBD}


\Copyright{The author. Released under the CC BY-ND license (International 4.0).}

\title{Matroids with bases as minimal resolving sets of graphs}

\author{Usman Ali$^{*,}$\authornote{1}
\and
Iffat Fida Hussain\authornote{2}
}

\authortext{1}{Centre for Advanced Studies in Pure and Applied Mathematics Bahauddin Zakariya University, Multan, Pakistan
   ($^*$Corresponding Author. \email{uali@bzu.edu.pk}) }
\authortext{2}{Centre for Advanced Studies in Pure and Applied Mathematics, Bahauddin Zakariya University, Multan, Pakistan
 (\email{fidahussainiffat7@gmail.com}).}

\begin{document}

\maketitle

\begin{abstract}
  We define an independence system associated with simple graphs. We prove that the independence system is a matroid for certain families of graphs, including trees, with bases as minimal resolving sets. Consequently, the greedy algorithm on the matroid can be used to find the minimum-cost resolving set of weighted graphs, wherein the independent system is a matroid. We also characterize hyperplanes of the matroid for trees and prove that its dual matroid is loop-free. 
\end{abstract}
\section{Introduction}
Let \emph{G(V, E)} denote a simple graph, comprised of a vertex set \emph{V(G)} and an edge set {$E(G)\subseteq V \times V$. Let $R\subseteq V(G)$, $R$ is designated a resolving set \cite{ch} or locating set \cite{slater} if every vertex of $G$ is uniquely defined by its distances from the vertices of $R$. It is designated a distinguishing set \cite{rall} if each vertex can be uniquely identified by its set of neighbors in $R$. Additionally, it is recognized as a locating dominating set \cite{rall} or beacon set \cite{col} if none of the sets of neighbors in $R$ is empty (that is, if $R$ is also a dominating set). A determining set in the context of a graph $G$ refers to a set of vertices $D$, such that every automorphism of 
$G$ is distinctly identified by its effect on the vertices of $D$. The minimum size of a determining set serves as a gauge for graph symmetry, playing a crucial role in addressing problems related to graph automorphisms \cite{alb, alb1}. The notion of a determining set was first introduced in \cite{bou} and independently as a fixing set in \cite{erw}. 
\begin{definition}
\label{Bou} \cite{boutin}
A set $S$ of vertices in a graph $G$ is res-independent (det-independent) if for every $ s \in S, S-\{s\} $ is not a resolving (determining) set.
\end{definition}
The above definition gives birth to an independence system. The author claims that with this definition, a maximal res-independent set is a minimal resolving (determining) set and the system becomes a matroid if exchange property holds for minimal resolving (determining) sets in a graph $G$. We give counter examples to the claim that every maximal independent set is minimal resolving (determining) set in \cite{Zill}. The independence system defined in \cite{boutin} is not a matroid even there is exchange property in minimal resolving (determining) sets because there exist maximal independent sets which are not minimal resolving (determining) sets. The novelty of the author's concept in \cite{boutin} of the independence system is evident by the fact that the primary motivation is to identify graphs for which the computation of the metric dimension or the determining number of a graph can be achieved within polynomial-time complexity. Notably, these problems are known to be NP-hard for general graphs \cite{balaji}. It is worth mentioning that \cite{boutin} introduced a fresh perspectives and idea that should not be overshadowed by a small oversight. This misstep, though important to address (see examples 1 and 2 in \cite{Zill}) does not diminish the overall impact and significance of her valuable contribution. The refined independence system defined in this paper, particularly when it manifests as a matroid for a given graph, presents a significant advancement. The matroid structure introduces an algorithm offering a more efficient computational approach to address the challenging problems of the metric dimension and determining number in comparison to the broader class of general graphs.
Our refined definition is as follow. 
 \begin{definition}
A set $S$ of vertices in a graph $G$ is independent if it is a subset of a minimal resolving (determining) set.
\label{def1}
\end{definition}
With this definition, two independence systems $(G)_{\text{res}}$ and $(G)_{\text{det}}$  are obtained with respect to minimal resolving sets and minimal determining sets of a graph respectively. In our current study, we confine our focus to $(G)_{\text{res}}=(V(G), I)$ where $V(G)$ is our ground set and $I$ is the collection of independent sets (Definition \ref{def1}). For the remainder of the article, our investigation centers on the study of $(G)_{\text{res}}$. The rank of the independence system $(G)_{\text{res}}$ is the metric dimension of the graph $G$ denoted by $\beta(G)$. \\
\indent While Khuller \cite{balaji} has previously established a polynomial-time complexity of computing metric dimension of trees, we contribute a new proof within the framework of matroid theory. Specifically, our result demonstrates that the metric dimension of some families of graphs including trees can be computed using the greedy algorithm based on a matroid. This alternative proof offers a fresh perspective within the context of matroid theory, emphasizing the applicability of efficient algorithms in solving the metric dimension problems for graphs whose independent systems are matroids.
 It is noteworthy that our proofs of matroids focus on demonstrating the augmentation property, a more suitable approach than relying on the base exchange property. Unlike the latter, which necessitates prior knowledge that all minimal resolving sets must share the same cardinality. Our preference for highlighting the augmentation property stems from its practicality, especially in cases where finding the uniform cardinality of minimal resolving sets for general graphs might be challenging. It is essential to note that while having the same cardinality is a necessary condition but not sufficient, as illustrated in the examples of wheel graph $W_{5}$ in Theorem \ref{wheel}(c) and $W_{7}$ in Theorem \ref{wheel}(e), it alone does not ensure that the independence system qualifies as a matroid.\\
    \indent The weighted version of the metric dimension of a graph was introduced in \cite{epstein} and gives algorithms for some graph to find the minimum-cost weighted minimum resolving set. Our independence system $(G)_{\text{res}}$ (when it forms a matroid) provides a more established frame work in the context of matroid theory to find the minimum-cost weighted minimum resolving set. More precisely, if $(G)_{\text{res}}$ is a matroid then the well-known greedy algorithm based on the matroid can be used to find the minimum-cost weighted base of the matroid. The greedy algorithm 1 do the job of finding minimum-cost. The set
 $\{1, 2,...,n\}$ in the algorithm represents the vertex set $\{ v_1,v_2,...,v_n\}$ of the graph $G$ when $(G)_{\text{res}}$ is a matroid.
 \newline
\begin{table}[h!]
   \centering
    \resizebox{\textwidth}{!}{%
\begin{tabular}{l}

\hline
\textbf{Algorithm 1 }Greedy algorithm for selecting the min-weight base of a
matroid \\ \hline
\textbf{Input:} the matroid $\mathbf{M}=\left( V(G),\mathbf{I}\right) ,$ where $%
V(G)=\left\{ 1,2,...,n\right\} $ is the ground set, and weight of $i$ is $w_{i}.
$ \\ 
\textbf{Output:} A base $B\in \mathbf{I}$ such that $w\left( B\right)
=\min_{B\in \beta }w\left( B\right) .$ \\ 
$1:$ Relabel the elements of the matroid so that $w_{1}\leq w_{2}\leq
...\leq w_{n}.$ \\ 
$2:$ $S\longleftarrow \emptyset .$ \\ 
$3:$ \textbf{for} $i\longleftarrow 1$ to $n$ \textbf{do} \\ 
$4:$ \ \ \ \textbf{if} $S+i\in \mathbf{I}$ \textbf{then} \\ 
$5:$ \ \ \ \ \ \ $S\longleftarrow S+i.$ \\ 
$6:$ \ \ \ \textbf{end if} \\ 
$7:$ \textbf{end for} \\ 
$8:$ \textbf{return} $S$ \\ \hline
\end{tabular}
}
\end{table}
\newline
\indent The presence of a loop in a matroid introduces challenges for determining crucial parameters, such as the chromatic number, which cannot be defined for a matroid with loops \cite{mc1, mc2}. A similar notion to an independence system in topology and combinatorial commutative algebra exists, known as an abstract simplicial complex \cite{stanly}. A loop-free independence system is exactly the same as the simplicial complex \cite{stanly}. When a simplicial complex adheres to the matroid structure, it is referred to as a matroid-complex and is therefore shellable. A shellable simplicial complex has the desirable property that the associated Stanley-Reisner ring is sequentially Cohen Macaulay. The independence system $(G)_{\text{res}}$ provides direction to the researchers in commutative algebra to study properties of associated rings.\\
\indent Beyond our emphasis on trees $T_{n}$ and the characterization of their hyperplanes in the matroid $(T_{n})_{\text{res}}$, our paper delves into the dual matroid of $(T_{n})_{\text{res}}$, showing its loop-free. Additionally, we discuss $(G)_{\text{res}}$ for wheel graphs, complete graphs, cycle graphs, and prove whether the corresponding independence systems are matroid or otherwise. Finally, we contribute to the research landscape by presenting a list of questions at the conclusion of our paper, inviting further exploration and investigation in the field.
\section{Preliminaries}
\subsection{Metric Dimension}
In a connected graph, there exists a path consisting of one or more
edges between every pair of vertices. The shortest possible path from a vertex \emph{u}
to another vertex \emph{v} is determined by the number of edges between $u$ and $v$ and its length is the distance between \emph{u} and
\emph{v}, denoted by $d(u, v)$. In a graph \emph{G}, the
\emph{metric code} of \emph{v} in \emph{V(G)} with respect to an
ordered set $ R = \{ h_1, h_2, h_3, \dots, h_m\} \subset V(G)$, is a
\emph{m}-tuple given below:
$$ c_v(R) = ( d(v, h_1), d(v, h_2), d(v, h_3), \dots, d(v, h_m)),$$
where \emph{R} is known as resolving set if for $u, v \in V$, $ u
\neq v \Rightarrow c_u(R) \neq c_v(R)$. A \emph{ metric basis} for \emph{G}
is a resolving set of the minimum cardinality in \emph{G} and the
number of elements of basis is called \emph{G$'$s} \emph{metric
dimension} denoted by \emph{$\beta(G)$} \cite{8}. For some variants of the metric dimension of different families of graphs, see \cite{e1, e2, e3, e4, e5}. The conceptualization of a resolving set and the metric dimension was
commenced by Slater \cite{slater}, and independently by Harary and Melter \cite{14}. The relevance and implementation of resolving sets and
the metric dimension emanates in various fields of study including broadcasting \cite{6}, chemical
bonding, robot navigation \cite{balaji}(where resolving sets are called landmark set), and geographical routing
protocols.\\
 
\subsection{Basics of Matroids}
An independence system $(E,I)$  is also referred to as a \textit{hereditary system} \cite{west}.
The set $ X $ is designated as an \textit{independent set} if $X\subset I$ and
classified as a \textit{dependent set} otherwise. The empty set $\phi $ is
inherently independent and, the set $E$ is considered dependent by definition. The concept of independence varies across different contexts, leading to diverse types of independence systems. In linear algebra, the
independence corresponds to the usual linear independence \cite{whitney}.
In the context of graph theory, for a simple undirected graph $G$, the property is
edge-independent, meaning a set of edges is independent if its induced
graph is acyclic \cite{west}. A matroid provides a generalization of linear independence in linear algebra. The formal definition of a matroid is presented below, which will be employed in subsequent discussions.
\begin{definition} [\cite{west}]
\label{matroid}
 An independence system $(E,I)$ comprising a family $I$
with subsets of a finite set $E$ (ground set) is recognized  as a \textit{matroid } denoted by $M= (E,I)$ if it
adheres to one of the following equivalent properties:\\
1. (Base exchange property) for any two maximal independent sets $B_1$ and
$B_2$ and for every $x\in B_{1}$, there exists $y\in B_{2}$ such
that $(B_{1}\setminus \{x\})\cup \{y\}$ is also a maximal
independent set.\\
2. (Augmentation property) for any two independent sets $A$ and $B$ with $\mid A \mid<\mid B\mid$, $A\cup \{x\}$ is independent for some $x\in B\setminus A$.
\end{definition}
A maximal independent within a matroid is referred to as \textit{base} of the matroid. Conversely, a minimal dependent set is known as a \textit{circuit} of the matroid. An individual element from $E$ that is dependent is specifically termed as a \textit{loop}. A matroid is classified as loop-free if it lacks any loop. An element from $E$ is identified as an \textit{isthmus} if it is part of every base of the matroid.
The maximum cardinality among independent sets within a matroid $M$ is denoted as the \textit{rank} of the matroid, symbolized as $r(M)$. 
The rank of a subset $X$ of $E$ is determined by the maximum cardinality among independent sets contained in $X$ and the overall rank $r(M)$ corresponds to the rank of $E$. In other words, $r:2^E\rightarrow \textsc{R}$ represents a submodular function. 
A \textit{flat} $F$ of a matroid $M$ is defined as a subset of $E$ such that $r(F)<r(F\cup \{x\})$ for all $x\in E\setminus F$ . By definition, $E$ itself is considered a flat of $M$. A uniform matroid denoted as $U_{n}^{l}=(S, I)$, is defined over a set $S$ with cardinality $n$. In this matroid, every subset of $S$ with a cardinality of at most $l$, is considered an independent set. Notably, each set comprising exactly $l$ elements serves as a base for $U_{n}^{l}$ and a circuit within this matroid consists of precisely $l+1$ elements \cite{oxley}. A \textit{hyperplane} in a matroid $M=(E, I)$ is a flat $H$ of $E$ for which $r(H)=r(M)-1$. For a matroid $M=(E, I)$, the dual matroid of $M$ denoted as $M^*=(E, I^*)$, is a matroid such that bases of $M^*$ are the complement in $E$ of the bases of $M$ \cite{oxley}. The circuits of a matroid need not to be necessarily of the same size unlike the bases of a matroid \cite{oxley}. Two matroids $M_1$ and $M_2$ are said to be isomorphic if there is a bijection between the elements of the ground set which maps independent sets to independent sets, (or equivalently circuits to circuits, or bases to bases) \cite{gordon, Recski}. It is important to maintain clarity by distinguishing between the finite set $E$ and the edge set $E(G)$ of a graph $G$.\\
\indent In this paper, $E$ specifically refers to the vertex set $V(G)$ of a graph $G$.
The independence systems discussed in the preceding subsection give rise to matroids. The initial system deals with linearly independent sets
 in a vector space and is called \textit{linear matroid}. The second system, where independent sets comprise acyclic sets of edges in a simple
 undirected graph, is known as \textit{graphic matroid} \cite{oxley}.
\section{Our Results}
\subsection{Trees}
To facilitate the reader's understanding of our new results in this section, we introduce relevant concepts from \cite{balaji} and \cite{slater} pertaining to trees. In the discussion, we focus on all trees except paths because for a path graph $P_{n}$ shown in Figure \ref{pathtree}(a), both $\{v_{1}\}, \{v_{2}, v_{3}\}$ are minimal resolving sets (maximal independent sets), $\forall n  \geq 4$. As maximal independent sets do not share the same cardinality, therefore, $(P_{n})_{res}$ is not a matroid $\forall n \geq 4$.
\begin{figure}
\begin{center}
\includegraphics [width=10cm]{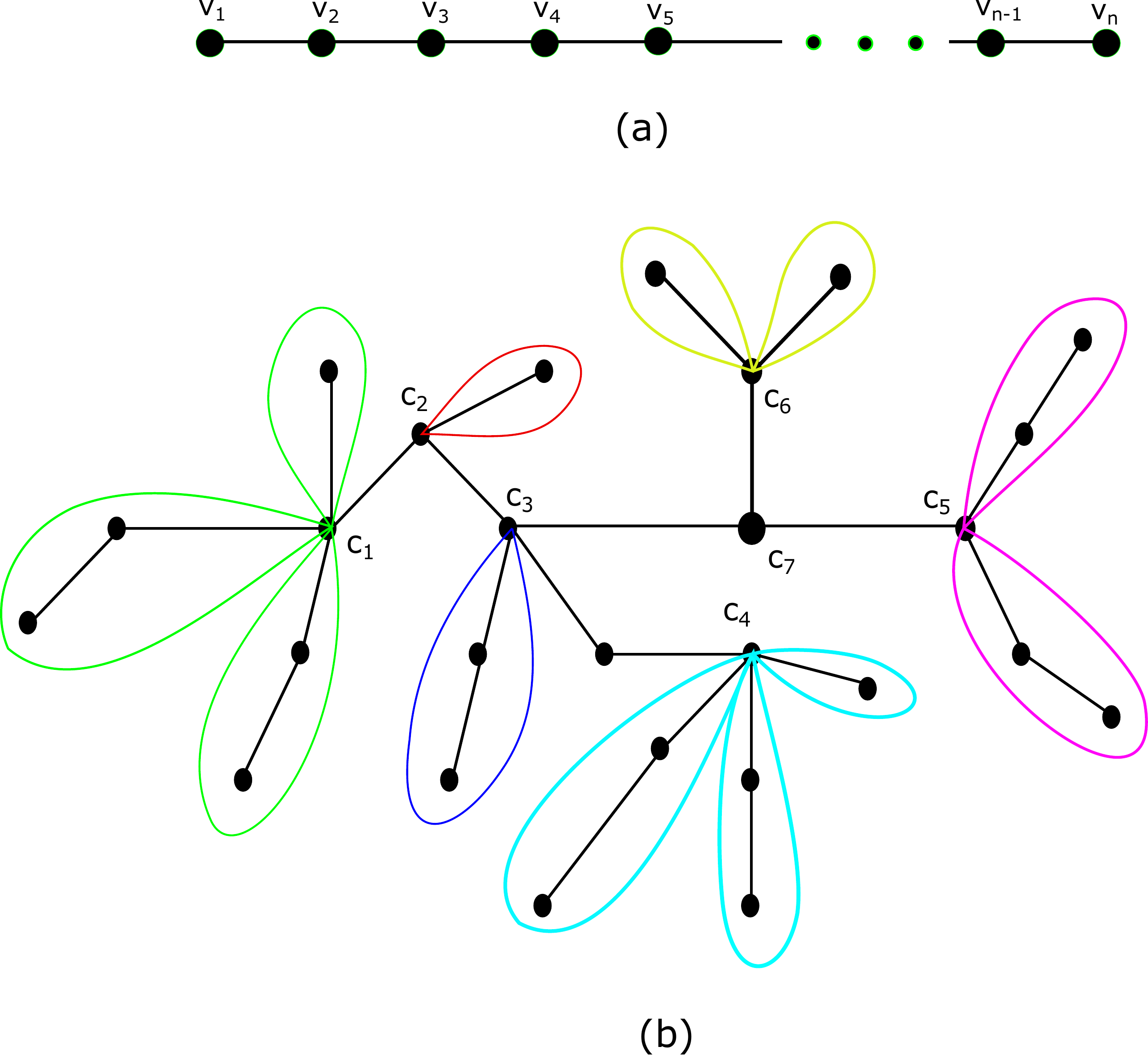}
\caption{(a) Path \qquad (b) Tree} \label{pathtree}
\end{center}
\end{figure} 
\begin{definition} [\cite{balaji, slater}]
In the context of a tree $T_n$ consisting of $n$ vertices, $\forall n \geq 4$, a branch at a specific vertex $v$ is defined as a connected component of $T_n-\{v\}\cup \{v\}$. Specifically, when $deg(v)\geq 3$, a branch that forms a path at $v$ is referred to as a branch (or a leg) of the tree at $v$.
\end{definition} 
\begin{definition} [\cite{balaji, slater}]
A vertex $v$ in a tree is designated as an exterior major vertex (or core vertex) if $deg(v) \geq 3$ and it possesses at least one branch path at $v$.
\end{definition}
In Figure \ref{pathtree}(b), a tree is depicted, featuring $6$ exterior major vertices. Each exterior major vertex has branch paths, which are highlighted with distinct colors. A vertex with a degree greater than or equal to three but lacking branch paths is termed a \textit{major not exterior} vertex. The distribution of vertices in the tree, categorized as exterior major and their associated branch paths, is presented in Table \ref{T1}.

 \begin{table}[h!]
\renewcommand\arraystretch{1.5}   
\noindent\[ 
\begin{array}{|c|c|c|c|}
\hline
\text{vertex} &\text{classification} & \text{number of branch paths} \\
\hline
c_1 & \text{exterior major} & 3 \\
\hline
c_{2} & \text{exterior major} & 1 \\
\hline
 c_{3}& \text{exterior major} & 1 \\
  \hline
c_{4} & \text{exterior major} & 3 \\ 
\hline
 c_{5} & \text{exterior major} & 2 \\
\hline
 c_{6} & \text{exterior major} & 2 \\ 
 \hline
 c_{7} & \text{major not exterior} & 0 \\
 \hline
\end{array}
\]
\caption{Vertex classification of the tree in Figure \ref{pathtree}(b).}\label{T1}
\end{table}
\begin{theorem} [\cite{slater}] (Minimal resolving set criterion).
\label{thm3}
A set $S$ of vertices is considered a minimal resolving set for a tree if and only if, for each exterior major vertex $v$, there exists a distinct vertex in $S$ from exactly all but one of the branch paths at $v$, and this vertex is different from $v$.
\end{theorem}
\begin{lemma}
  For any tree $T_{n}$ of order $n\geq 4$, it is guaranteed to have an exterior major vertex with atleast two branch paths.
  \label{l1}
\end{lemma}
\begin{proof}
Suppose the contrary, assuming that there is no exterior major vertex with two branch paths in tree $T_n$, and every exterior major vertex has only one branch path. Let $S$ be a minimal resolving set of $T_n$ then by Theorem \ref{thm3}, when we exclude one branch path of each exterior major from $T_n$, we are left with $S= \{\}$ which is impossible as $\beta(G) \geq 1$ and it is known that the metric dimension of a graph is $1$ if and only if it is path \cite{comp}. 
\end{proof}
\begin{theorem} 
\label{tree}
The independence system $(T_{n})_{\text{res}}$ forms a matroid for each tree $T_{n},  \forall \, n \geq 4$.
\end{theorem}
\begin{proof}
Consider two sets $X$ and $Y$, assumed to be independent with $|X|<|Y|$. It is important to note that $X$ cannot serve as a minimal resolving set because every minimal resolving set is a maximal independent set. Since $|X|<|Y|$, therefore, there exists $y\in Y$ and $y\notin X$ where $y$ is a vertex of a branch path of an exterior major vertex $v$ having atleast two branch paths (such an exterior major vertex is guaranteed by Lemma\ref{l1}). By Theorem \ref{thm3}, $X \cup \{y\}$ must necessarily be a subset of a minimal resolving set. The augmentation property of Definition 
\ref{matroid} is satisfied and hence, $(T_{n})_{\text{res}}$ is a matroid. 
\end{proof}

We characterize the hyperplanes of the matroid  $(T_{n})_{\text{res}}$ in the following theorem.
\begin{theorem} 
\label{hyperplane}
A subset $H$ of $V(T_n)$ is a hyperplane of $(T_{n})_{\text{res}}$ if and only if $V(H) = V(T_n) \setminus (X_{1}\cup X_{2})$, where $X_{1}, X_{2}$ are the vertex sets of two distinct branch paths of the same exterior major vertex of $T_n, \forall\, n\geq 4$.
\label{hyper}
\end{theorem}
\begin{proof}
First, suppose that $H$ is a hyperplane of $(T_{n})_{\text{res}}$. Consequently, $H$ is a flat and $r(H)= \beta (T_n)-1$ as $r((T_{n})_{\text{res}})= \beta (T_n)$. We suppose the contrary case where $X_{1}, X_{2}$ are  vertex sets of two different branch paths of the same exterior major of $T_n$ and $V(H)$ has non-empty intersection with $X_{1} \cup X_{2}$. In other words, $V(H)$ contains one vertex either from $X_{1}$ or $X_{2}$. According to Theorem \ref{thm3}, $H$ contains a minimal resolving set and therefore, $r(H)=\beta(T_n)$ leading to a contradiction. Thus, the only possibility is $V(H) = V(T_n) \setminus (X_{1}\cup X_{2})$.\\
Conversely, Suppose $V(H) = V(T_n) \setminus (X_{1}\cup X_{2})$. Again, utilizing  Theorem \ref{thm3}, if we augment $V(H)$ by some $ x \in X_{1}\cup X_{2}$, then $V(H) \cup \{x\}$ contains a minimal resolving set, and we obtain $r(H \cup \{x\} )= \beta(T_n)$. Consequently, $r(H \cup \{x\}) > r(H), \forall \, x \in X_{1}\cup X_{2}$, establishing that $H$ is a flat and $r(H)= \beta (T)-1$. Consequently, $H$ indeed becomes a hyperplane.
\end{proof}
\begin{figure}
\begin{center}
\includegraphics [width=15cm]{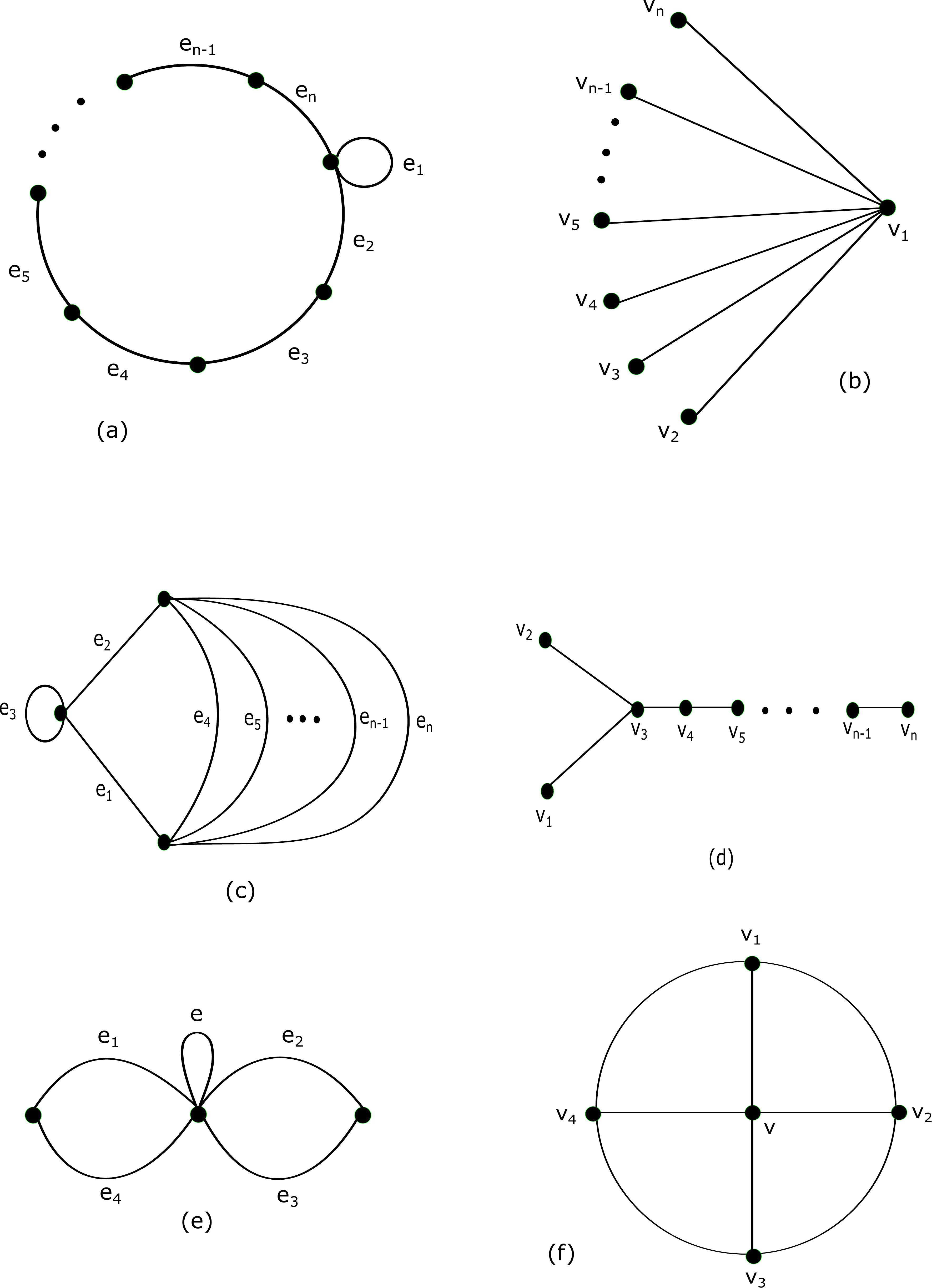}
\caption{}\label{fig1}
\end{center}
\end{figure}
\begin{proposition} 
\begin{itemize}
\qquad \qquad \qquad\qquad\qquad\qquad \qquad \qquad\qquad\qquad\qquad\qquad\qquad\qquad\qquad\qquad\qquad\qquad \\
\textbf{(a)} The matroid $(T_{n})_{\text{res}}$ has always a loop $\forall \, n \geq 4$.\\
\textbf{(b)} The matroid $(T_{n})_{\text{res}}$ has no isthmus $\forall \, n \geq 4$.
\end{itemize}
\label{propo1}
\end{proposition}
\begin{proof}
\textbf{(a)}. According to Lemma \ref{l1}, every tree invariably possesses an exterior major vertex (say $v$) with two branch paths. This exterior major vertex $v$ does not fall in any minimal resolving set, as per the criterion outlined in Theorem \ref{thm3}. This is because each minimal resolving set includes a vertex from a branch path of $v$ but not $v$ itself. Consequently, there always exists a singleton set $\{v\}$ that is not a subset of any minimal resolving set. Hence, $(T_{n})_{\text{res}}$ consistently harbors a singleton dependent set or a loop.\\
\textbf{(b)} As per the criterion outlined in Theorem \ref{thm3}, every minimal resolving set of a tree $T_n$ incorporates a vertex from all branch paths (excluding one) of each exterior major of $T$. Consequently, there is no single vertex that must be included necessarily in each minimal resolving set. This implies there is no vertex that always belongs to every base of $(T_{n})_{\text{res}}$. As a result, $(T_{n})_{\text{res}}$ lacks any isthmus.  
\end{proof}
For a matroid $M= (E, I)$, $E\backslash H$ is recognized a circuit of the dual Matroid $M^*$ if and only if $H$ serves as a hyperplane in Matroid $M$, as detailed in Proposition 3.18 on page 116 of \cite{gordon}. The subsequent theorem elucidates that the dual matroid $T^*_{\text{res}}$  is devoid of any loop. 
\begin{theorem} 
\label{dual}
The dual matroid $(T_{n})^*_{\text{res}}$ is loop-free, $\forall \, n \geq 4$.
\end{theorem}
\begin{proof}
By the definition of dual matroid, the bases of $(T_{n})^*_{\text{res}}$ are the complements of minimal resolving sets of $T_n$. We know that $V(T_n) \setminus H$ is a circuit of $(T_{n})^*_{\text{res}}$ if and only if $H$ is a hyperplane in $(T_{n})_{\text{res}}$.  According to Theorem \ref{hyper}, each hyperplane  $H$ in $(T_{n})_{\text{res}}$ possesses the property that its vertex set $V(H)$ is given by $V(T_n) \backslash (X_{1}\cup X_{2})$, where $X_{1}, X_{2}$ represent the vertex sets of two different branch paths originating from the same exterior major of $T_n$. Consequently, every circuit $C=X_{1}\cup X_{2}$ in $(T_{n})^*_{\text{res}}$ must have the cardinality $\geq 1$, precluding the existence of singleton dependent sets. As a result, it follows that $(T_{n})^*_{\text{res}}$ is loop-free.  
\end{proof}
\begin{proposition}
\qquad \qquad \qquad\qquad\qquad\qquad \qquad \qquad\qquad\qquad\qquad\qquad\qquad\qquad\qquad\qquad\qquad\qquad \\ 
\label{cycle matroid}
\textbf{({i})} The matroid $(T_{n})_{\text{res}}$ associated with the tree depicted in Figure \ref{fig1}(b) is isomorphic to the graphic matroid corresponding to the connected graph illustrated in Figure \ref{fig1}(a).\\
\textbf{({ii})} The matroid $(T_{n})_{\text{res}}$ associated with the tree depicted in Figure \ref{fig1}(d) is isomorphic to the graphic matroid corresponding to the connected graph illustrated in Figure \ref{fig1}(c).
\end{proposition}
\begin{proof}
\textbf{(\emph{i})}. For the tree $T_n$ shown in Figure \ref{fig1}(b), its matroid $(T_{n})_{\text{res}}= (V(T_n), I)$, where $I$ contains independent sets (subsets of minimal resolving sets of the tree). In the case of the graph $G$ depicted in Figure \ref{fig1}(a), the graphic matroid is denoted as $M=(E_{1}, I_{1})$, where $ E_{1}=\{e_{1}, e_{2}, e_{3},\ldots, e_{n}\}$ and $I_{1}$ represents the family of acyclic edges in the graph $G$. The bijective map $\phi : E_{1}\rightarrow V(T_n) $
$$\phi(e_{i})=v_{i}, \qquad \text{for all} \qquad  1\leq i \leq n. $$  maps the independent sets of $M$ to the independent sets of $(T_{n})_{\text{res}}$. Consequently, $\phi$ establishes an isomorphism between $M$ and $T_{\text{res}}$.\\
\textbf{(\emph{ii})}. The same bijective map $\phi$ defined for part \textbf{(\emph{i})} establishes an isomorphism between the two matroids.
\end{proof}

\subsection{Wheels}
The wheel graph on $n+1$ vertices, denoted as $W_{n}$, is formed by taking the join of the cycle $C_{n}$ and the complete graph $K_{1}$ for $n \geq 3$. In other words, $W_{n}$ consists of an $n$-cycle with an additional vertex $v$ that is adjacent to all the vertices of the cycle shown in Figure \ref{figure5}. Specifically, the metric dimensions of $W_{n}$ for $n \geq 3$ are listed below in Table \ref{TW} \cite{8}.
\begin{figure}
\begin{center}
\includegraphics [width=6cm]{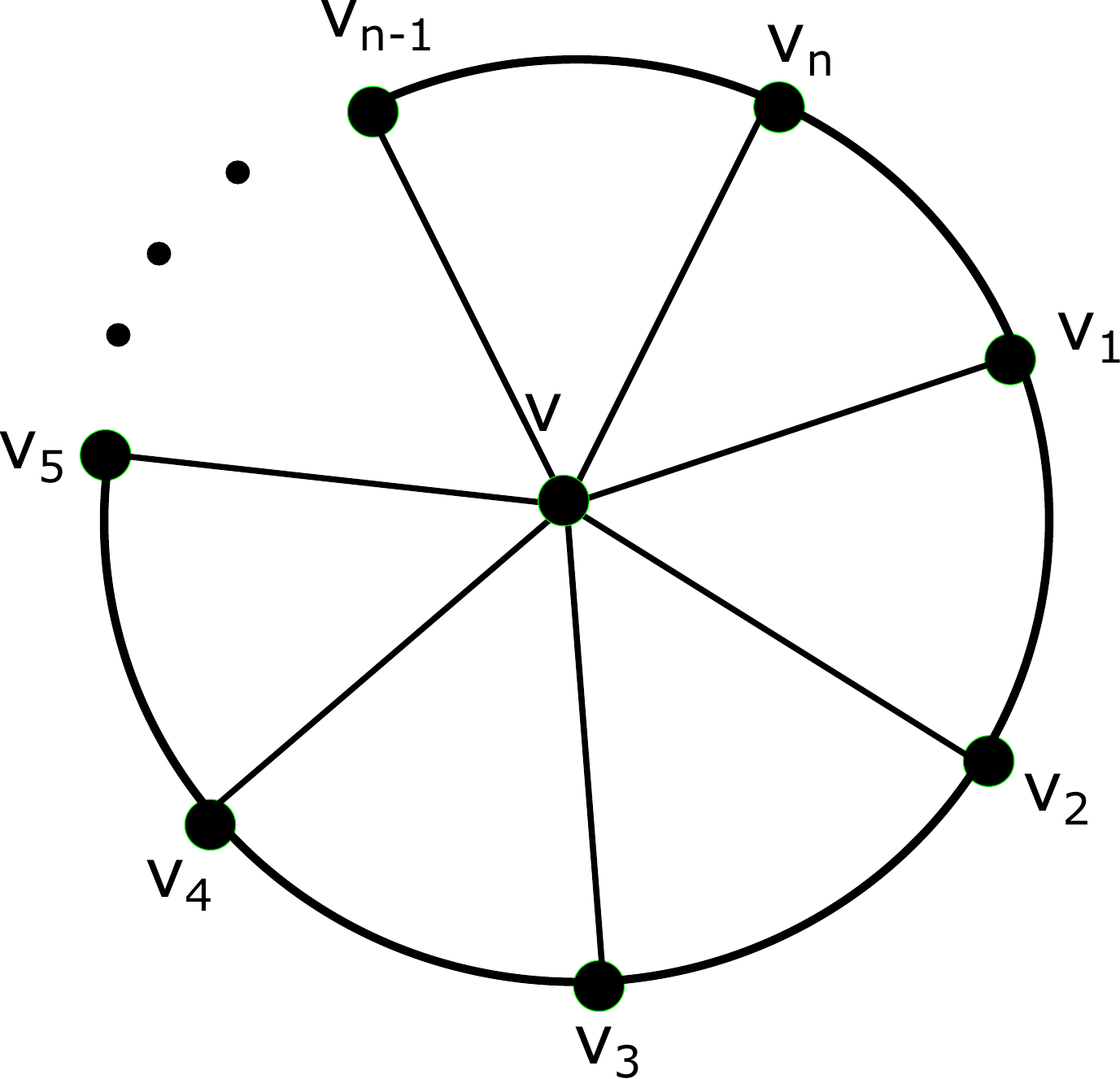}
\caption{Wheel $W_{n}$}
 \label{figure5}
\end{center}
\end{figure}

\begin{table}[h!]
\begin{tabular}[h]{|c| c| }

\hline
$W_{n}$& $\beta (W_{n})$ \\
\hline
$W_3$ & $3$ \\
$W_{4}$ & $2$\\
 $W_{5}$& $2$ \\
$W_{6}$ & $3$ \\ 
 $W_{n}, \forall \, n \geq 7$ & $\lfloor \frac{2n+2}{5} \rfloor$ \\
  \hline
\end{tabular}
\centering\caption{Metric dimension table of $W_{n}$}
\label{TW}
\end{table}
\begin{table}[h!]

\begin{tabular}[h]{|c| c| c|c| c| c| }

\hline
vertex& $v_{1}$ & $ v_{2}$& $ v_{3}$& $ v_{4}$& $ v$ \\
\hline
$v_1$ & 0 & 1&2&1&1 \\
$v_{2}$ & 1 & 0&1&2&1\\
 $v_{3}$& 2 & 1&0&1&1 \\
$v_{4}$ & 1 & 2&1&0&1 \\ 
 $v$ & 1 & 1&1&1&0 \\
  \hline
\end{tabular}
\centering\caption{Distance table of $W_{4}$}
\label{T2}
\end{table}
\begin{table}[h]

\begin{tabular}[h]{|c| c| c|c| c| c|c| }

\hline
vertex& $v_{1}$ & $ v_{2}$& $ v_{3}$& $ v_{4}$& $ v_{5}$ &$v$ \\
\hline
$v_1$ &   0 & 1 &2&2&1&1    \\
$v_{2}$ & 1 & 0&1&2&2&1    \\
 $v_{3}$& 2 & 1&0&1&2&1    \\
$v_{4}$ & 2 & 2&1&0&1&1    \\ 
 $v_5$ & 1 & 2&2&1&0&1     \\
 $v$&1&1&1&1&1&0           \\
  \hline
\end{tabular}
\centering\caption{Distance table of $W_{5}$}
\label{T3}
\end{table}

\begin{table}[h!]

\begin{tabular}[h]{|c| c| c|c| c| c|c|c| }

\hline
vertex& $v_{1}$ & $ v_{2}$& $ v_{3}$& $ v_{4}$& $ v_{5}$ &$v_{6}$ &$v$ \\
\hline
$v_1$ &   0 & 1  & 2 & 2 &2 & 1 & 1    \\
$v_{2}$& 1 & 0 &  1 & 2 &2 & 2 &1    \\
 $v_{3}$& 2 & 1&0&1&2&2&1    \\
$v_{4}$ & 2 & 2&1&0&1&2&1    \\ 
 $v_5$ & 2 & 2&2&1&0&1&1     \\
 $v_6$&1&2&2&2&1&0&1           \\
 $v$&1&1&1&1&1&1&0\\
  \hline
\end{tabular}
\centering\caption{Distance table of $W_{6}$}
\label{T4}
\end{table}
\begin{table}[h!]

\begin{tabular}[h]{|c| c| c|c| c| c|c|c|c| }

\hline
vertex& $v_{1}$ & $ v_{2}$& $ v_{3}$& $ v_{4}$& $ v_{5}$ &$v_{6}$ &$v_7$&$v$ \\
\hline
$v_1$ &   0 & 1  & 2 & 2 & 2 &2 & 1 & 1    \\
$v_{2}$& 1 & 0 &  1 & 2 &2&2 & 2 &1    \\
 $v_{3}$& 2 & 1&0&1&2&2&2&1    \\
$v_{4}$ & 2 & 2&1&0&1&2&2&1    \\ 
 $v_5$ & 2 &2& 2&2&1&0&1&1     \\
 $v_6$&1&2&2&2&2&1&0&1           \\
 $v_7$&1&2&2&2&2&1&0&1\\
 $v$&1&1&1&1&1&1&1&0\\
  \hline
\end{tabular}
\centering\caption{Distance table of $W_{7}$}
\label{T5}
\end{table}
\newpage
\begin{theorem} 
\qquad \qquad \qquad\qquad\qquad\qquad \qquad \qquad\qquad\qquad\qquad\qquad\qquad\qquad\qquad\qquad\qquad\qquad \\ 
  (a) $(W_{3})_{res}$ is a matroid. \\
  (b) $(W_{4})_{res}$ is a matroid.\\
  (c) $(W_{5})_{res}$ is not a matroid.\\
  (d) $(W_{6})_{res}$ is a matroid.\\
  (e) $(W_{7})_{res}$ is not a matroid.\\
  (f) $(W_{n})_{res}, \forall \, n \geq 8$ is not a matroid.
  \label{wheel}
\end{theorem}
\begin{proof}
(a) Since $W_{3}= K_{4}$ and every $3$ cardinality subset of $V(W_{3})$ is a minimal resolving set. Consequently, $(W_{3})_{res}$ is a uniform matroid of rank $3$. \\
(b) For wheel $W_{4}$, Table \ref{T2} lists all distances.The minimal resolving sets for the graph $W_{4}$ are $\{v_{1}, v_{2}\}$, $\{v_{1}, v_{4}\}$, $\{v_{2}, v_{3}\}$, and $\{v_{3}, v_{4}\}$. The independence system of $(W_{4})$ is:
 $(W_{4})_{res}=\{ \phi, \{v_{1}\}, \{v_{2}\}, \{v_{3}\}, \{v_{4}\}, \{v_{1}, v_{2}\}, \{v_{1}, v_{4}\}, \{v_{2}, v_{3}\}, \{v_{3}, v_{4}\}$.
  It is noteworthy that the independent sets exhibit the augmentation property. Consequently, $(W_{4})_{res}$ is established as a matroid. Importantly, $(W_{4})_{res}$ is characterized as a non-uniform matroid, given the presence of two cardinality sets of vertices $\{v_{1}, v_{3}\},\{v_{2}, v_{4}\}$ that do not qualify as independent sets. Moreover, the matroid $(W_{4})_{res}$ of the wheel $W_{4}$ shown in Figure \ref{fig1}(f) is isomorphic to the graphic matroid of the connected graph shown in Figure \ref{fig1}(e). \\
 (c) For wheel $W_{5}$, Table \ref{T3} lists all distances. The minimal resolving sets for the graph $W_{5}$ are composed of the pairs $\{v_{1}, v_{2}\}$, $\{v_{1}, v_{5}\}$, $\{v_{2}, v_{3}\}$, and $\{v_{3}, v_{4}\}$ and $\{v_{4}, v_{5}\}$. It is noteworthy that the independent sets lack the augmentation property. For example, for the sets $\{v_{1}\}$ and $\{v_{3}, v_{4}\}$, by augmentation we get $\{v_{1}, v_{3}\}$ and  $\{v_{1}, v_{4}\}$ which do not qualify as independent sets. Consequently, $(W_{5})_{res}$ is not a matroid. It is one of the rare cases where maximal independent sets share the same cardinality but still the independent system is not a matroid.\\
 (d) For wheel $W_{6}$, Table \ref{T4} lists all distances.Importantly, $(W_{6})_{res}$ is characterized as a uniform matroid of rank $3$ as all sets of cardinality at most three qualify as independent sets.\\
 (e) For graph $W_{7}$, Table \ref{T5} lists all distances. The minimal resolving sets for $W_{7}$ are: $\{v_{1}, v_{2}, v_{4}\},\{v_{1}, v_{2}, v_{6}\}, \{v_{1}, v_{3}, v_{4}\}, \{v_{1}, v_{3}, v_{5}\},
  \{v_{1}, v_{3}, v_{6}\}, \{v_{1}, v_{3}, v_{7}\}, \{v_{1}, v_{4}, v_{6}\}, \{v_{1}, v_{5}, v_{6}\},\\
  \{v_{1}, v_{5}, v_{7}\}, \{v_{2}, v_{4}, v_{5}\}, \{v_{2}, v_{4}, v_{6}\}, \{v_{2}, v_{4}, v_{7}\}, \{v_{2}, v_{5}, v_{7}\}, \{v_{2}, v_{6}, v_{7}\}, \{v_{3}, v_{4}, v_{6}\} , \{v_{3}, v_{5}, v_{6}\},\\ \{v_{3}, v_{5}, v_{7}\}, \{v_{4}, v_{5}, v_{7}\}, \{v_{4}, v_{6}, v_{7}\}$. \\It is noteworthy that the independent sets lack augmentation property. For example, if we consider two independent sets $\{v_{1}, v_{2}\}$ and $\{v_{3}, v_{5}, v_{7}\}$, upon augmenting $\{v_{1}, v_{2}\}$ by one of the vertices $v_{3}, v_{5}$, and $v_{7}$, we get $\{v_{1}, v_{2}, v_{3}\}, \{v_{1}, v_{2}, v_{5}\}, \{v_{1}, v_{2}, v_{7}\}$ respectively,
  which are not subsets of any minimal resolving set.  Consequently, $(W_{7})_{res}$ is not a matroid although all minimal resolving sets share the same cardinality.\\
 (e) Boutin proved in Theorem $7$ of \cite{boutin} that $\forall \, n \geq 8$, minimal resolving sets in $W_{n}$ do not possess exchange property. Consequently, $(W_{n})_{res}, \forall \, n \geq 8$ is not a matroid.
 \end{proof}
\subsection{Cycles and Complete Graphs} 
The cycle graph $C_{n}$, comprising $n$ vertices is characterized by a single cycle traversing all vertices. Given that each vertex in $C_{n}$ has a degree of $2$, it follows that no individual vertex can resolve all the vertices. The minimal resolving set for $C_{n}$ is of cardinality $2$, and $\beta(G)=2$ according to \cite{comp}. Notably, for even values of $n$, the cycle graph $C_{n}$ transforms into an antipodal graph, where each vertex has an antipode. An antipode of a vertex $v$ is defined as a vertex $u$ situated at the maximum distance from $v$ within the graph. The set of two antipodal vertices is not a resolving set for $C_{n}$.  
\begin{theorem}
The independent system $(C_{n})_{\text{res}}$ is a matroid for a cycle graph $C_{n}, \forall \, n \geq 3$.\\
\end{theorem}
  \begin{proof}
  \textbf{Case 1}: When $n$ is odd \\ 
For the set $V(C_{n})$ where $|V(C_{n})|=n$, every subset of cardinality $2$ serves as a minimal resolving set for $C_{n}$. According to Definition \ref{def1}, every set with a cardinality less than or equal to $2$ qualifies as an independent set for $(C_n)_{\text{res}}$. In other words, all $2$ cardinality sets of vertices are maximal independent sets (i.e., bases). As a result, $(C_{n})_{\text{res}}$ is characterized as a uniform matroid of rank $2$. \\

\textbf{Case 2}: When $n$ is even \\

Consider sets $X$ and $Y$ assumed to be independent, with $|X|<|Y|$. This implies that $Y$ must have cardinality $2$.  Consequently, $Y$ is a minimal resolving set and the following two subcases arise:\\

\textbf{Case (a)}: $X\subset Y$ \\

If $X$ is augmented by adding a vertex $v$ from $Y$ for $v \in Y \setminus X$, then $X \cup \{v\}$ becomes a minimal resolving set. The augmentation property in Definition \ref{matroid} is satisfied and hence $(C_{n})_{\text{res}}$ is a matroid of rank $2$. \\

\textbf{Case (b):} $X \nsubseteq Y$ \\
If $X$ is augmented by a vertex $v$ from $Y$ such that $v$ is not antipodal to the element in $X$ (such a vertex must exist because of the uniqueness of antipodal vertex), then $X \cup \{v\}$ is necessarily a minimal resolving set.
Moreover, if ${v_{1}}$ and ${v_{2}}$ represent two antipodal vertices in $C_{n}$, then the set $\{v_{1}, v_{2}\}$ does not constitute a minimal resolving set and is not a subset of any other minimal resolving set. This implies that $(C_{n})_{\text{res}}$ is a non-uniform matroid of rank $2$.
\end{proof}

Every vertex in complete graph $K_{n}$ is adjacent to every other vertex. As a result, each vertex can uniquely identify only itself and is at the same distance from all other vertices. To discern all vertices in the graph, a resolving set with the minimum cardinality of $n-1$ is essential. Hence, $\beta(K_{n}) = n-1$ \cite{comp}, and any set of vertices with a cardinality of $n-1$ serves as a minimal resolving set. According to Definition \ref{def1}, every set of cardinality less than or equal to $n-1$ is an independent set in $(K_{n})_{\text{res}}$. In other words, all sets of cardinality $n-1$ are maximal independent sets (i.e., bases). Consequently, $(K_{n})_{\text{res}}$ is a uniform matroid of rank $n-1$.\\
A bipartite graph $G$ with $n$ vertices is defined by partitioning its vertex set into two sets, $U$ and $W$, such that $|U|=m$, $|W|=n$. Importantly, each edge $v_{1}v_{2}$ in $G$ is characterized by the condition that $v_{1}$ and $v_{2}$ cannot belong to the same subset, either $U$ or $W$ within the vertex set $V(G)$. A bipartite graph that is complete, denoted as $K_{m, n}$, is one where all elements in set $U$ share the same neighborhoods, and the same holds true for all elements in set $W$. For a minimal resolving set in $K_{m, n}$, it is essential to include elements from both sets $U$ and $W$ symmetrically (excluding exactly one element from each set). As a consequence, the cardinality of a minimal resolving set must be $m+n-2$, and $\beta(K_{m, n})=m+n-2$ \cite{comp}.
 Consequently, $(K_{m, n})_{res}$ is a matroid of rank $m+n-2$.
\section{Open Questions}
In Proposition \ref{cycle matroid}, we provide two examples: one illustrating the isomorphism of the graphic matroid of the graph in Figure \ref{fig1}(a) and the matroid $(T_n)_{\text{res}}$ of the given tree in Figure \ref{fig1}(b); and the second example features the tree shown in Figure \ref{fig1}(d), whose matroid $(T_n)_{\text{res}}$ is isomorphic to the graphic matroid of the connected graph shown in Figure \ref{fig1}(c). This leads to the formulation of an open question:

\textbf{Question 1}. Characterize all trees for which $(T_n)_{\text{res}}$ is isomorphic to a graphic matroid of some connected graph.

The metric dimension of a Wheel graph can be solved by a linear-time algorithm \cite{epstein}. We prove that $(W_{n})_{\text{res}}$, for all $n \geq 7$, is not a matroid.

\textbf{Question 2}. Does there exist other families of graphs for which the metric dimension problem is polynomial-time solvable, but $(G)_{\text{res}}$ is not a matroid?

\textbf{Question 3}. Investigate whether $(G)_{\text{res}}$ and $(G)_{\text{det}}$ for different families of graphs are matroids or otherwise.

\textbf{Question 4}. Study properties of rings associated with $(G)_{\text{res}}$ and $(G)_{\text{det}}$.


\begin{thebibliography}{999}
\bibitem{alb} M. O. Albertson, D. L. Boutin. Using determining sets to distinguish Kneser graphs. \emph{Electron. J. Combin.} 14: 1-9, 2007.
\bibitem{alb1} M. O. Albertson, D. L. Boutin. Automorphisms and distinguishing numbers of geometric cliques. \emph{Disc. Comput. Geom.} 39: 778-785, 2008.
\bibitem{6} Z. Beerloiva, F. Eberhard, T. Erlebach, A. Hall, M. Hoffmann. Network discovery and verification. \emph{IEEE J. Sel. Areas Commun.} 24: 2168-2181, 2006.  
\bibitem{bou}  D. L. Boutin. Identifying graph automorphisms using determining sets. \emph{Electron. J. Combin.} 13: 1-12, 2006.
\bibitem{boutin} D. L. Boutin. Determining Sets, Resolving Sets, and the Exchange Property. \emph{Graphs Comb.} 25: 789-806, 2009. 
\bibitem{8} P. S. Buczkowski, G. Chartrand, C. Poisson, P. Zhang. On k-dimensional graphs and their basis. \emph{Period. Math. Hung.} 46: 9-15, 2003. 
\bibitem{ch} G. Chartrand, L. Eroh, M. A. Johnson, O. R. Oellermann. Resolvability in graphs and the metric dimension of a graph. \emph{Discret. Appl. Math.} 105: 99-113, 2000. 
\bibitem{comp} G. Chartrand, L. Eroh, M. A.Johnson, O. R. Oellermann. Resolvability in graphs and the metric dimension of a graph. \emph{Discret. Appl. Math.} 105: 99-113, 2000. 
\bibitem{col} C. I. Colbourn, P. J. Slater, L. K. Stewart. Locating dominating sets in series parallel networks. \emph{Congr. Numer.} 56: 135-162, 1987.
\bibitem{erw} D. J. Erwin, F. Harary. Destroying automorphisms by fixing nodes. \emph{Discret. Math.} 306: 3244-3252, 2006.
\bibitem{epstein} L. Epstein, A. Levin, G. J. Woeginger. The (Weighted) Metric Dimension of Graphs: Hard and Easy Cases. \emph{Algorithmica}. 72: 1130–1171, 2015.
\bibitem{e1} M. Fazil, I. Javaid, M. Salman, U. Ali. Locating-dominating sets in hypergraphs. \emph{Period. Math. Hung.} 72: 224-234, 2016. 
\bibitem{gordon} G. Gordon, J. McNulty.  Matroids: A geometric introduction. Cambridge University Press, 2012.
\bibitem{14} F. Harary, R. A. Melter. On the metric dimension of a graph. \emph{Ars Combin.} 2: 191-195, 1976.
\bibitem{e2} I. F. Hussain, S. Afridi, A. M. Qureshi, G. Ali, U. Ali. Fault-tolerant resolvability of swapped optical transpose
 interconnection system. \emph{J. Math.} 2022: 1-6, 2022.
\bibitem{e3} I. Javaid, M. Fazil, U. Ali, M. Salman. On some parameters related to fixing sets in graphs. \emph{J. Prim. Res. Math.} 14: 1-12, 2018.
\bibitem{balaji} S. Khuller, P. Raghavachari, A. Rosenfeld. Landmarks in graphs. \emph{Discret. Appl. Math.} 70: 217-229, 1996. 
\bibitem{mc1} M. Lason. Indicated coloring of matroids. \emph{Discret. Appl. Math.} 179: 241-243, 2014.
\bibitem{mc2} A. Lehman. A solution to Shannons switching game. \emph{SIAM J}. 12: 687-725, 1964.
\bibitem{e4} S. Naz, M. Salman, U. Ali, I. Javaid, S. A. Bokhary. On the constant metric dimension of generalized peterson graph $P(n, 4)$. Acta. Math. Sin.-English Ser. 30: 1145–1160, (2014).
\bibitem{oxley} J. G. Oxley. Matroid Theory. Oxford University Press, Oxford, 1992. 
\bibitem{rall} D. F. Rall, P. J. Slater.  On location-domination number for certain classes of graphs. \emph{Congr. Numr.} 45: 97-106, 1984.
\bibitem{Recski} A. Recski. Matroid theory and applications: in electric networks and in statics. Springer-Verlag Berlin Heidelberg. 1989.
\bibitem{Zill} Z. E. Shams, M. Salman, U. Ali.  On Maximal Det-Independent(Res-Independent) Sets in Graphs. \emph{Graphs Comb.} 2022: 1-5, 2022. 
\bibitem{slater} P. J. Slater. Leaves of trees. \emph{Congr. Numer.} 14: 549-559, 1975.
\bibitem{stanly} R. P. Stanley. Combinatorics and Commutative algebra. Birkhauser, 1983. 
\bibitem{e5} R. C. Tillquist, R. M. Frongillo, M. E. Liadser. Getting the lay of the land in discrete space: a survey of metric dimension and its applications. \emph{SIAM Review} 65: 919-962, 2023.  
\bibitem{west} D. B. West. Introduction to Graph Theory. Prentice Hall, 2001.
\bibitem{whitney} H. Whitney. On the abstract properties of linear dependence. \emph{Amer. J. Math.} 57: 509-533, 1935. 
\end{thebibliography}
\end{document}